\documentclass[11pt]{article}
\usepackage[centertags]{amsmath}
\usepackage{amsfonts}
\usepackage{amssymb}
\usepackage{amsthm}
\usepackage{newlfont}
\usepackage{graphicx}
\usepackage{multirow}
\usepackage{rotating}

\theoremstyle{definition}
\newtheorem{definition}{Definition}
\theoremstyle{remark}

\theoremstyle{theorem}
\newtheorem{theorem}{Theorem}

\newtheorem{lemma}{Lemma}
\begin{document}
\thispagestyle{empty}

%
%

\title{ Control of accuracy on Taylor-collocation method to solve the weakly regular Volterra
integral equations of the first kind by using the CESTAC method}
\author{ Samad Noeiaghdam, $^{a,}$
\footnote{E-mail addresses: s.noeiaghdam.sci@iauctb.ac.ir;
samadnoeiaghdam@gmail.com}~~~~Denis Sidorov $^{b,c,}$
\footnote{E-mail addresses: contact.dns@gmail.com} ~~and~~Valery
Sizikov $^{d,}$ \footnote{E-mail addresses: sizikov2000@mail.ru} }
  \date{}
 \maketitle

\begin{center}
\scriptsize{ $^{a}$Department of Mathematics, Central Tehran Branch, Islamic Azad University, Tehran, Iran. \\
$^{b}$Irkutsk National Research Technical University, Irkutsk, Russia.\\
 $^{c}$Energy Systems Institute of Russian Academy of Sciences.\\
 $^{d}$ ITMO University, St.Petersburg, Russia.\\
}
\end{center}
\begin{abstract}
Finding the optimal parameters and functions of iterative methods is
among the main problems of the Numerical Analysis. For this aim, a
technique of the stochastic arithmetic (SA) is used to control of
accuracy on Taylor-collocation method for solving first kind weakly
regular integral equations (IEs). Thus, the CESTAC
\footnote{Controle et Estimation Stochastique des Arrondis de
Calculs} method is applied and instead of usual mathematical
softwares the CADNA \footnote{Control of Accuracy and Debugging for
Numerical Applications} library is used. Also, the convergence
theorem of presented method is illustrated. In order to apply the
CESTAC method we will prove a theorem that it will be our licence to
use the new termination criterion instead of traditional absolute
error. By using this theorem we can show that number of common
significant digits (NCSDs) between two successive approximations are
almost equal to NCSDs between exact and numerical solution. Finally,
some examples are solved by using the Taylor-collocation method
based on the CESTAC method. Several tables of numerical solutions
based on the both arithmetics are presented. Comparison between
number of iterations are demonstrated by using the floating point
arithmetic (FPA) for different values of $\varepsilon$.

\vspace{.5cm}{\it Keywords: Stochastic arithmetic, Floating point
arithmetic, CESTAC method, CADNA library,
 Taylor-collocation method, First kind Volterra integral equation.}
\end{abstract}
\section{Introduction}
The first kind IEs  are among of the important problems in the
Applied Mathematics which have many applications in mathematical
sciences and engineering
\cite{man1,man2,man5,man6,dn1,man7,1,dn2,man9,2,siz1,siz2,siz3,siz4,siz5}.
Many authors have investigated analytical or numerical methods to
solve the first kind IEs
\cite{6a,man6,4,3,man7,dn4,man9,siz3,siz4,siz5}. Since finding the
solution of the first kind IEs by using analytical methods is
challenging, especially for applications, thus the numerical methods
suggested to solve the approximate solutions specially for ill-posed
or singular IEs of the first kind.

In order to solve the the first kind IEs we have several numerical
methods such as collocation method \cite{man1,5,man9,2}, Adomian
decomposition method \cite{8,7,6}, homotopy analysis method
\cite{man2,man5,man6,man10,man7} and many others
\cite{13,10,11,9,12,siz1,siz2,siz3,siz4,siz5}.

The aim of this work is to investigate the following  first kind
classical Volterra IE
\begin{equation}\label{1}
\int_{0}^{t} k(r,s) v(s) ds  = f(r),~ 0 \leq r \leq T \leq 1,
\end{equation}
where the kernel $k(r,s)$ is discontinuous along continuous curves
$\rho_i, i=1,2,\cdots,m-1$ as
\begin{equation}\label{2}\left\{
\begin{array}{l}
  k_1(r,s),~~~~~ 0= \rho_0(r)  \leq s \leq  \rho_1(r),\\
    \\
  k_2(r,s),~~~~~ \rho_1(r)  \leq s \leq  \rho_2(r),\\
    ~~~\vdots\\
    k_m(r,s),~~~~~ \rho_{m-1}(r)  \leq s \leq  \rho_m(r) =r \leq 1. \\
\end{array}\right.
\end{equation}
Finally, we can rewrite Eq. (\ref{1}) with conditions (\ref{2}) as
follows
\begin{equation}\label{3}
\int_{\rho_0(r)}^{\rho_1(r)} k_1(r,s) v(s) ds +
\int_{\rho_1(r)}^{\rho_2(r)} k_2(r,s) v(s) ds + \cdots +
\int_{\rho_{m-1}(r)}^{\rho_m(r)} k_m(r,s) v(s) ds = f(r),
\end{equation}
where $f(0)=0$. The pressed form of Eq. (\ref{3}) is obtained in the
following form
\begin{equation}\label{4}
\sum_{p=1}^m \int_{\rho_{p-1}(r)}^{\rho_p(r)} k_p(r,s) v(s) ds =
f(r).
\end{equation}

In the last decades, many mathematical schemes have been presented
to solve the IEs with jump discontinuous kernels
\cite{dav,dn8,dn1,dn5,dn6,dn3,dn4,dn7,dn2}. In this paper, by
combining the collocation method with Taylor polynomials the
Volterra IE (\ref{4}) is solved. Recently, many linear and
non-linear problems have solved by using Taylor-collocation method
\cite{16,15,18,14,19,17}.

We know that the mentioned researches are accomplished by using the
FPA. In these papers, the numerical results were presented for a
special iteration. In \cite{siz6} (see also [58, p. 274]), when
solving IEs of the type $Av = f$ by the iteration method: $v_n = (I-
A^{\ast}A) v_{n-1} + A^{\ast}f, n = 1, 2, \cdots$, the following
rules for stopping the iteration process were proposed:

1) by the discrepancy, namely, stop at such a number $n$, for which
$\|Av_n - f \| \leq \delta + \xi \|v_n\|$, where
$\|\tilde{f}-f\|\leq \delta, \|\tilde{A}-A\|\leq \xi$ ($\delta$ and
$\xi$ are the errors of $f$ and $A$);

2) by the correction, namely, $\|v_n-v_{n-1}\| \leq a_1 \delta + a_2
\xi$ ($a_1> 0$ and $a_2> 0$ are some numbers). It is proved that
$\|v-v_n\| \rightarrow 0$ for $\delta, \xi \rightarrow 0$; however,
for finite $\delta, \xi$, finding $\varepsilon$ is difficult.
\\
Also, generally the efficiency of numerical methods were
investigated by applying the traditional absolute error as follows
\begin{equation}\label{cond1}
\left| v(s) - v_n(s) \right| \leq \varepsilon.
\end{equation}
which depends on the exact solution, the small tolerance value
$\varepsilon$ and the number of iterations. Now, we can ask how
authors can find that which iteration is suitable to stop?  If we do
not know the exact solution, how we can apply the condition
(\ref{cond1})? How we can find the proper value of $\varepsilon$? In
this condition, for large values of $\varepsilon$ we can not find
the appropriate approximation. Because the iterations will be
stopped without getting to suitable solution. Also, for small values
of $\varepsilon$ we will have useless iterations without amending
the precision. Thus, according to mentioned problems we can not
accept the condition (\ref{cond1}).

In this paper, the new arithmetic is suggested instead of
traditional arithmetic which is called the SA
\cite{cad7,cad8,cad9,cad10,cad11,cad12}. Thus, we will apply the
Taylor-collocation method to solve the first kind Volterra IE
(\ref{4}) based on the SA. Also, the validate and control of
accuracy on numerical results are investigated by using the CESTAC
method \cite{cad1,cad2,cad3,cad4,cad5,cad6} and the CADNA library
\cite{cad15,cad14,cad13}. By using this method, we can find the
optimal approximation, the optimal iteration, the optimal error of
numerical methods \cite{man3,man4,man8,man10}. Also, the CESTAC
method does not have the disadvantages of the FPA. The stopping
condition of the CESTAC method is not similar to the FPA. It depends
on two successive approximations instead of usual absolute error and
in order to apply the new condition we will prove a theorem which
will show the NCSDs between exact and approximate solutions are
almost equal to NCSDs between two successive approximations. In new
arithmetic, we are applying the CADNA library instead of
mathematical softwares like Maple, Matlab, Mathematica and so on
\cite{cad15,cad14, man4,cad13,man8,man10}. The CADNA library is a
programming environment that it has many abilities. The
 CADNA programs can be written by using C, C++ or FORTRAN
 codes \cite{man3,man4,man8,man10}. Also, some of numerical instabilities such as  bifurcations, blow-up, mathematical operations, branching,
functions and others can be found by using the CADNA library.

This paper is organized in the following form: Section 2 presents
the Taylor-collocation method and employ it for solution of  the
first kind IE (\ref{4}). Also, the convergence theorem of presented
method is proved. Section 3 describes the application of the SA. The
CESTAC method and the CADNA library are considered  in full details.
Furthermore, a sample program of the CADNA library is presented  for
methodological objectives. Numerical validation of
Taylor-collocation method based on the CESTAC method is investigated
in Section 4. In this section, the new applicable termination
criteria for the SA is presented. The theorem which is proved in
this section, is our license to use the new stopping condition. In
this theorem, we will show the NCSDs between two iterations are
almost equal to the NCSDs between exact and approximate solutions.
The CADNA algorithm to solve the IE. (\ref{4}) is presented in
Section 5. Also, several examples are solved based on the presented
algorithm. The numerical results are obtained for both environments,
the FPA and the SA. In this section, the optimal iteration,
approximation and error of Taylor-collocation method are presented.
Finally, the conclusions and some advantages of the  SA than the FPA
are investigated in Section 6.

\section{Taylor-collocation method}
Let  us search for the solution in the following form
\begin{equation}\label{5}
 v_n(s) = \sum_{j=0}^n \frac{1}{j!} v^{(j)} (c) (s-c)^j + O(h^{n+1}),
\end{equation}
which is the Taylor polynomial of degree $n$ at point $s=c$ where
\begin{equation}\label{2222}
e_n(s) = |v(s)-v_n(s)| = O(h^{n+1}).
\end{equation}
Substituting Eq. (\ref{5}) into Eq. (\ref{4}) leads to
\begin{equation}\label{6}
\sum_{p=1}^m \int_{\rho_{p-1}(r)}^{\rho_p(r)} k_p(r,s) \sum_{j=0}^n
\frac{1}{j!} v^{(j)} (c) (s-c)^j ds = f(r).
\end{equation}
By producing the collocation points
\begin{equation}\label{2-1}
r_i=a+(\frac{b-a}{n})i, ~~~~~i=0,1,\cdots, n,
\end{equation}
and putting grids (\ref{2-1}) in Eq. (\ref{6}) we get
\begin{equation}\label{2-2}
\sum_{j=0}^n \frac{1}{j!} \left[ \sum_{p=1}^m
\int_{\rho_{p-1}(r_i)}^{\rho_p(r_i)} k_p(r_i,s)
  (s-c)^j ds \right] v^{(j)} (c) = f(r_i).
\end{equation}

Now, we can write Eq. (\ref{2-2}) in the following form
\begin{equation}\label{2-3}
AV=F,
\end{equation}
where
$$
A=\left[
\begin{array}{cccc}
  A_{00} & A_{01} & \cdots & A_{0n} \\
  A_{10} & A_{11} & \cdots & A_{1n} \\
  \vdots & \vdots & \ddots & \vdots \\
  A_{n0} & A_{n1} & \cdots & A_{nn}
\end{array}
\right]_{(n+1)(n+1)},
$$
that
$$
A_{ij} = \sum_{p=1}^m \int_{\rho_{p-1}(r_i)}^{\rho_p(r_i)}
k_p(r_i,s)
  (s-c)^j ds,
$$
and
$$
V=\left[ \begin{array}{cccc}
           \bar{v}^{(0)}(c) & \bar{v}^{(1)}(c) & \cdots & \bar{v}^{(n)}(c)
         \end{array}
   \right]^T,
$$
$$
F=\left[ \begin{array}{cccc}
           F(r_0) & F(r_1)  & \cdots & F(r_n)
         \end{array}
   \right]^T.
$$
Coefficients $\bar{v}^{(j)}(c)$ are uniquely determined by system
(\ref{2-3}). So Eq. (\ref{4}) has only one unique solution which can
be obtained by
\begin{equation}\label{8}
 \bar{v}_n(s) = \sum_{j=0}^n \frac{1}{j!} \bar{v}^{(j)} (c) (s-c)^j.
\end{equation}

\subsection{Convergence analysis}

The main theorem of convergence for presented method to find the
approximate solution of weakly regular Volterra IE (\ref{4}) is
presented. At first, several lemmas are required that we should
remaind. Also, we note that $I = [0,1]$, $X = L^2(I)$.

\begin{lemma}\label{L1}\cite{20}
Let $X = L^2(I)$ and $W$ be a Volterra integral operator on $X$ with
square summable kernel $k(r,s)$ where $\int_0^1\int_0^1 |k(r,s)| dr
ds = M^2$, and $M$ is a constant. Then the operator $W$ can be
defined as
$$
W(\phi(r)) = \int_0^r k(r,s) \phi(s) ds,
$$
which $W$ is bounded and
$$
\|W(\phi(r)) \| \leq M \| \phi \|.
$$
\end{lemma}

\begin{lemma}\label{L2} \cite{21} Let $v(s)$ be a sufficiently smooth function on $[0, 1]$ and
$p_n(s)$ is the interpolating polynomial to $v(s)$ at points $s_i$
then
$$
v(s) - p_n(s) = \frac{v^{n+1}(s)}{(n+1)!} \prod_{i=0}^n (s-s_i),
s\in[0,1],
$$
where $s_i, i = 0, 1, . . . , n$ are the roots of interpolating
polynomial. We can write
$$
| v(s) - p_n(s) | \leq \frac{M_n}{2^{2n+1} (n+1)! },
$$
where
$$
M_n = \max \left\{ |v^{n+1}(s)|; s \in (0,1) \right \}.
$$
\end{lemma}

\begin{lemma}\label{L3}
The linear operators $Z_j, j=1,2,...,m$ can be defined on $X$ as
$$
\begin{array}{l}
\displaystyle  Z_1(v(r)) = \int_{\rho_0(r)}^{\rho_1(r)} k_1(r,s)v(s) ds,  \\
    \\
\displaystyle   Z_2(v(r)) = \int_{\rho_1(r)}^{\rho_2(r)} k_2(r,s)v(s) ds,  \\
 ~~~~~   \vdots\\
\displaystyle     Z_m(v(r)) = \int_{\rho_{m-1}(r)}^{\rho_m(r)} k_m(r,s)v(s) ds,  \\
\end{array}
$$
where $s \in [0,1]$ and $v \in L^2(I)$.
\end{lemma}

According to Lemma \ref{L1}, the mentioned operators $Z_1, Z_2,
\cdots, Z_m$ are bounded. For one-one and onto $Z_1$ the inverse
operator $Z_1^{-1}$ is bounded.

\begin{theorem} Let
$ v_n(s) = \sum_{m=0}^n \frac{1}{m!} v^{(m)} (c) (s-c)^m$ be the
approximate solution of  Eq. (\ref{1}) which is obtained from
applying the Taylor-collocation method and $ \bar{v}_n(s) =
\sum_{m=0}^n \frac{1}{m!} \bar{v}^{(m)} (c) (s-c)^m$ is the
expansion of exact solution $v(s)$. Thus
\begin{equation}\label{9}
  \|e_n(s)\|_{2} \leq (\mu_1 \mu_2+\cdots+\mu_1 \mu_m)\left(\frac{M_n}{2^{2n+1} (n+1)! }+\gamma_n
 \|Q_n\|_2\right),
\end{equation}
where $M_n = \max \left\{ |v^{n+1}(s)|; s \in (0,1) \right \}$ and
$Q_n = (q_0(c),q_1(c),\cdots,q_n(c)); q_n = v^{(n)}(c) -
\bar{v}^{(n)}(c)$.
\end{theorem}

\textbf{Proof:} By substituting presented operators of Lemma
\ref{L3} in Eq. (\ref{4}) we have
\begin{equation}\label{4-1}
Z_1(v(r)) +  Z_2(v(r)) + \cdots + Z_m(v(r)) = f(r),
\end{equation}
and by using inverse operator $Z_1^{-1}$ we get
\begin{equation}\label{4-2}
v(r) + Z_1^{-1}[ Z_2(v(r))] + \cdots + Z_1^{-1}[Z_m(v(r))] =
Z_1^{-1} [f(r)].
\end{equation}
Also, Eq. (\ref{4-2}) for $n$-th order approximation can be written
in the following form
\begin{equation}\label{4-3}
v_n(r) + Z_1^{-1}[ Z_2(v_n(r))] + \cdots + Z_1^{-1}[Z_m(v_n(r))] =
Z_1^{-1} [f(r)].
\end{equation}
From subtracting Eqs. (\ref{4-2}) and (\ref{4-3}) we have
\begin{equation}\label{4-4}
e_n(r) = Z_1^{-1}[ Z_2(v_n(r)-v(r))] + \cdots +
Z_1^{-1}[Z_m(v_n(r)-v(r))].
\end{equation}
By applying Lemma \ref{L3} and using Eq. (\ref{4-4}), there are
constant values $\mu_1,\cdots, \mu_m$ such that
\begin{equation}\label{4-5}
\begin{array}{ll}
  e_n(r)& \leq \mu_1 \mu_2\|v(r) - v_n(r)\|_2 + \cdots + \mu_1 \mu_m \|v(r) -
v_n(r)\|_2 \\
   \\
    &\leq (\mu_1\mu_2+\cdots+\mu_1\mu_m) \|v(r) - v_n(r)\|_2.
\end{array}
\end{equation}

Now, for exact and approximate solutions  $v(s)$ and $v_n(s)$ we can
write
\begin{equation}\label{10}
\| v(s) - v_n(s) \|_{2} =  \| v(s)- \bar{v}_n(s)+\bar{v}_n(s)-
v_n(s) \|_{2} \leq \| v(s) - \bar{v}_n(s) \|_{2} +\| \bar{v}_n(s) -
v_n(s) \|_{2},
\end{equation}
where
$$
v_n(s) = \sum_{m=0}^n \frac{1}{m!} v^{(m)} (c) (s-c)^m,
$$
and
$$
\bar{v}_n(s) = \sum_{m=0}^n \frac{1}{m!} \bar{v}^{(m)} (c) (s-c)^m.
$$

From Lemma \ref{L2} we get
\begin{equation}\label{4-6}
\| v(s) - \bar{v}_n(s) \|_{2} \leq \frac{M_n}{2^{2n+1} (n+1)! },
\end{equation}
where
$$
M_n = \max \left\{ |v^{n+1}(s)|; s \in (0,1) \right \},
$$
and
\begin{equation}\label{4-7}
\| \bar{v}_n(s) - v_n(s) \|_{2} = \left[ \int_0^1 \left(
\sum_{m=0}^n (v^{(n)}(c) - \bar{v}^{(n)}(c) ) \frac{(s-c)^m}{m!}
\right)^2 ds \right]^{\frac{1}{2}} \leq \gamma_n \| Q_n\|_2,
\end{equation}
where
$$
\gamma_n = \left(  \sum_{m=0}^n \int_0^1  \left| \frac{(s-c)^m}{m!}
\right|^2 ds \right)^{\frac{1}{2}},
$$
and
$$Q_n = (q_0(c),q_1(c),\cdots,q_n(c)); q_n = v^{(n)}(c) -
\bar{v}^{(n)}(c).$$

Thus from Eqs. (\ref{10}), (\ref{4-6}) and (\ref{4-7}) we have
\begin{equation}\label{4-8}
\| v(s) - v_n(s) \|_{2} \leq \frac{M_n}{2^{2n+1} (n+1)! } + \gamma_n
\| Q_n\|_2.
\end{equation}
Finally, by using Eqs. (\ref{4-5}) and (\ref{4-8}) we get
$$
 \|e_n(s)\|_{2} \leq (\mu_1 \mu_2+\cdots+\mu_1 \mu_m)\left(\frac{M_n}{2^{2n+1} (n+1)! }+\gamma_n
 \|Q_n\|_2\right).
$$
 Theorem is proved.


\section{Stochastic arithmetic}

The CESTAC method and the CADNA library are important method and
tool to simulate the Taylor-collocation method for solving problem
(\ref{1}). In this section, the main part of CESTAC method and the
CADNA library are investigated. More definitions, properties,
abilities and applications of the CESTAC method and the CADNA
library can be found in \cite{cad7,cad8,cad9,cad10,cad11,cad12}.

\subsection{The CESTAC methodology}
Assume that $F$ is the set of representable values which are
generated by computer. So for $g \in \mathbb{R}$ there is $G \in F$
which $G$ shows the reproduced form of $g$ by computer. Now, for $P$
mantissa bits of the binary FPA we have
\begin{equation}\label{0019}
G=g-\chi 2^{E-P}\gamma,
\end{equation}
where the sign of $g$ showed by $\chi$, the lost part of the
mantissa due to round-off error demonstrated by $2^{-P}\gamma$ and
the binary exponent of the result displayed by $E$. Also, we note
that in computer systems for single and double precisions we have
 $P = 24, 53$ respectively \cite{cad1,cad2,cad3,cad4,cad5,cad6}.

If the value $\gamma$ considers as a stochastic variable uniformly
distributed on $[-1, 1]$ then we can make perturbation on last
mantissa bit of $g$. Thus the obtained results for $G$ will be a
random variable with mean $(\mu)$ and the standard deviation
$(\sigma)$. The accuracy of random variable $G$ depends on both
parameters $(\mu)$ and $(\sigma)$ \cite{cad10,cad11,cad12}.

By $l$ times performing the process for $G_i , i = 1, . . . , l$ the
distribution of them is in the quasi Gaussian form. Therefore, the
mean of them is equal with the exact value of $g$ and the values of
$\mu$ and $\sigma$ can be estimated by these $l$ samples . The
following algorithm of the CESTAC method is presented where
$\tau_{\delta}$ is the value of $T$ distribution with $l-1$ degree
of freedom and confidence interval $1-\delta$.

\textbf{Algorithm 1:}

\texttt{Step 1- Find $l$ samples for $G$ as $ \Phi = \left\{ G_1,
G_2, ..., G_l \right\}$ by means of the perturbation of the last bit
of mantissa.}

\texttt{Step 2- Compute $\displaystyle G_{ave} = \frac{\sum_{i=1}^l
G_i}{l}$.}

\texttt{Step 3- Calculate $\displaystyle \sigma^2=
\frac{\sum_{i=1}^l ( G_i - G_{ave})^2}{l-1}$. }

\texttt{Step 4- Find the NCSDs between $G$ and $G_{ave}$ by using
$\displaystyle C_{G_{ave}, G} = \log_{10}\frac{\sqrt{l}
\left|G_{ave} \right|}{\tau_{\delta} \sigma}$.}

\texttt{Step 5- If $C_{G_{ave}, G} \leq 0$ or $G_{ave}=0,$
then write $G=@.0$.}\\

\subsection{The CADNA library}

Since the mathematical packages such as Mathematica, Maple and
others do not have this ability to produce the random variables, so
we should  introduce the new arithmetic and the novel software for
stochastic computations. To this aim, the CADNA library is
introduced. By using this software the SA can be applied instead of
the FPA. The CADNA library should run on LINUX operating system and
the commands of CADNA library must be written by using C, C++,
FORTRAN or ADA codes \cite{man3,man4,man8,man10}.

By using the CADNA library we can control the round-off error
propagation. Also, detecting the numerical instabilities is one of
the  important applications of this package. The CADNA algorithm
will be stopped when the NCSDs of numerical results equals to zero
and it will be shown by informatical zero $@.0$. Thus, we can find
the optimal iteration, approximation and error by using the
informatical zero \cite{cad1,cad2,cad3,cad4,cad5,cad6}.

The following program is a sample of CADNA library that we can find
some important remarks about abilities of this library.\\
program sample\\
\texttt{$\sharp$include <cadna.h>}\\
\texttt{cadna$_{-}$init(-1);\\}
\texttt{ main()\\}
\texttt{\{} \\
   \texttt{double$_{-}$st VALUE;} \\
\texttt{do}\\
    \texttt{\{}\\
The Main Program;\\
\texttt{printf(" \%s  ",strp(VALUE));}\\
\texttt{\}}\\
\texttt{ while(v[n]-v[n-1]!=0);\\}
 \texttt{cadna$_{-}$end();\\}
 \texttt{\}}

By applying \texttt{$\sharp$include <cadna.h>} we can call the
commands of CADNA library. Since all parameters are in the
stochastic mode, we should define them in \texttt{double},
\texttt{int} or other forms with \texttt{$_{-}$st} as
\texttt{double$_{-}$st}. Also, because the parameters defined based
on the SA thus in output we should apply \texttt{$\%s$} and
\texttt{strp} in print line. We note that in order to stop the
presented program, the new termination criterion is applied which
depends on difference of two successive approximations. Finally,
command \texttt{cadna$_{-}$end()} is the end of CADNA library.

\section{Numerical validation of Taylor-collocation method}

In order to solve the mentioned problems of the FPA, we suggest a
novel stopping condition based on the SA and the CESTAC method. This
condition is in the following form
\begin{equation}\label{cond2}
\left| v_{n+1} (s) -  v_n (s) \right| = @.0,
\end{equation}
which depends on two successive approximations $v_{n+1} (s)$ and
$v_n (s)$ and will be stopped when the difference of $v_{n+1} (s)$
and $v_n (s)$ equals to informatical zero $@.0$.

Now, what is our licence to apply the new termination criterion
(\ref{cond2}) instead of stopping condition (\ref{cond1})? In order
to apply the new condition, we should prove a new theorem to show
the equality  of the NCSDs between exact and approximate solutions
and the NCSDs between two successive approximations.

\begin{definition} \label{def2} \cite{cad1,cad2,cad3}
For numbers $\theta_1, \theta_2 \in \mathbb{R}$ the NCSDs can be
obtained as follows\\
 (1) for $\theta_1\neq \theta_2$,
\begin{equation}\label{13}
C_{\theta_1,\theta_2}=\log_{10}\left|\frac{\theta_1+\theta_2}{2(\theta_1-\theta_2)}
\right| = \log_{10}\left|\frac{\theta_1}{\theta_1-\theta_2} -
\frac{1}{2} \right|,
\end{equation}
(2) for all real numbers $\theta_1$, $C_{\theta_1,\theta_1} =
+\infty$.
\end{definition}

\begin{theorem} \label{th6}
Let $v(r)$ and $v_n(r)$ be the exact and approximate solutions of
weakly regular Volterra IE (\ref{4}) then
\begin{equation}\label{14}
\displaystyle C_{v_{n}(r),v(r)} - C_{v_{n}(r),v_{n+1}(r)} =
\mathcal{O}\left( h^{n+1}\right),
\end{equation}
where $C_{v_{n}(r),v(r)}$ shows the NCSDs of $v_{n}(r), v(r)$ and
$C_{v_{n}(r),v_{n+1}(r)}$ is the the NCSDs of two iterations
$v_{n}(r),v_{n+1}(r)$.
\end{theorem}

\textbf{Proof:} According to Definition \ref{def2}

$$
\begin{array}{l}
\displaystyle C_{v_n(r),v_{n+1}(r)} =\log_{10}\left|
\frac{v_n(r)+v_{n+1}(r)}{2(v_n(r)-v_{n+1}(r))}\right| =
\log_{10}\left|
\frac{v_n(r)}{v_n(r)-v_{n+1}(r)} - \frac{1}{2} \right| \\
\\
~~~~~~~~~~~~~~~~~~~~\displaystyle
=\log_{10}\left|\frac{v_n(r)}{v_n(r)-v_{n+1}(r)} \right| +
\log_{10}\left|1- \frac{1}{2
v_n(r)}(v_n(r)-v_{n+1}(r)) \right|\\
\\
~~~~~~~~~~~~~~~~~~~~\displaystyle
=\log_{10}\left|\frac{v_n(r)}{v_n(r)-v_{n+1}(r)} \right| +
\mathcal{O} \left(v_n(r)-v_{n+1}(r)\right).
\end{array}
$$

Since,
$$
v_n(r)-v_{n+1}(r) = v_n(r)-v(r)-(v_{n+1}(r)-v(r)) = E_n (r) -
E_{n+1}(r),
$$
thus,
$$
 \mathcal{O}\left(v_n(r)-v_{n+1}(r)\right) = \mathcal{O}\left(E_n (r) - E_{n+1}(r)\right) = \mathcal{O}\left( h^{n+1} \right) +  \mathcal{O}\left( h^{n+2} \right)= \mathcal{O}\left( h^{n+1} \right).
$$
Therefore,
\begin{equation}\label{15}
C_{v_n(r),v_{n+1}(r)} =
\log_{10}\left|\frac{v_n(r)}{v_n(r)-v_{n+1}(r)} \right|
+\mathcal{O}\left(h^{n+1}\right).
\end{equation}
Furthermore,
\begin{equation}\label{16}
 \begin{array}{ll}
\displaystyle C_{v_n(r),v(r)} &\displaystyle = \log_{10}\left|
\frac{v_n(r)+v(r)}{2(v_n(r)-v(r))}\right| = \log_{10}\left|
\frac{v_n(r)}{v_n(r)-v(r)} - \frac{1}{2} \right| \\
\\
&\displaystyle =\log_{10}\left| \frac{v_n(r)}{v_n(r)-v(r)}\right| +
\mathcal{O}(v_n(r) -
v(r))\\
\\
&\displaystyle = \log_{10}\left| \frac{v_n(r)}{v_n(r)-v(r)}\right|
+\mathcal{O}\left(h^{n+1}\right).
\end{array}
\end{equation}

By using Eqs. (\ref{15}) and (\ref{16}) we get,
$$\begin{array}{ll}
\displaystyle C_{v_n(r),v(r)} - C_{v_n(r),v_{n+1}(r)} &\displaystyle
= \log_{10}\left| \frac{v_n(r)}{v_n(r)-v(r)}\right| -
\log_{10}\left|\frac{v_n(r)}{v_n(r)-v_{n+1}(r)}
\right| +\mathcal{O}\left(h^{n+1}\right)\\
\\
&\displaystyle = \log_{10}\left| \frac{v_n(r)-v(r)
}{v_n(r)-v_{n+1}(r) } \right| +\mathcal{O}\left(h^{n+1}\right)\\
\\
&\displaystyle = \log_{10}\left|
\frac{\mathcal{O}\left(h^{n+1}\right)}{\mathcal{O}\left(h^{n+1}\right)} \right| +\mathcal{O}\left(h^{n+1}\right)\\
\\
&\displaystyle = \mathcal{O}\left(h^{n+1}\right),
  \end{array}
$$
and finally the following formula can be obtained
$$
\displaystyle C_{v_n(r),v(r)} - C_{v_n(r),v_{n+1}(r)} =
\mathcal{O}\left(h^{n+1}\right).
$$

When $n \rightarrow \infty$ the right hand side of above equation
tends to zero and hence we obtain
$$
C_{v_n(r),v(r)} = C_{v_n(r),v_{n+1}(r)}.
$$

Theorem \ref{th6} shows that the NCSDs between two approximations
$v_n(r),v_{n+1}(r)$ are almost equal to the NCSDs of $v_n(r),v(r)$.
Thus, we can apply termination criteria (\ref{cond2}) instead of
condition (\ref{cond1}).

\section{Numerical Illustration}
In this section, several examples of weakly regular Volterra IEs of
the first kind are presented. These examples are solved by using the
Taylor-collocation method based on the FPA and the SA. The CESTAC
method is applied to find the optimal approximation, the optimal
step and the optimal error of presented method. Also, both
conditions (\ref{cond1}) and (\ref{cond2}) are investigated to
compare the numerical results. It is clear that the termination
criterion (\ref{cond1}) depends on parameter $\varepsilon$ and
specially the exact solution and the number of iteration $n$. But
the algorithm of CESTAC method only depends on two successive
approximations. By comparison between these results, we will show
that the presented condition based on the SA is more applicable than
criterion (\ref{cond1}) which is based on the FPA. The following
CESTAC algorithm is presented to solve the problems.

\noindent \textbf{Algorithm 2:}\\
\texttt{Step 1- Let $n=1$.}\\
\texttt{Step 2- Define the Taylor polynomials and kernels $k_i(r,s)$. }\\
\texttt{Step 3- Enter the parameters $a, b, t$. }\\
\texttt{Step 4- Do the following steps while $|v_{n+1}(r) -
v_{n}(r)| \neq$ @.0 }

\texttt{\{}

\texttt{Step 4-1- Calculate the collocation points. }

\texttt{Step 4-2- Approximate the coefficient matrix $A$.}

\texttt{Step 4-3- Construct the right hand side matrix $F$.}

\texttt{Step 4-4- Solve the system (\ref{2-3}) by using the obtained
results of steps (4-2) and (4-3).}

\texttt{Step 4-5- Apply Eq. (\ref{5}) to find the numerical solution
of Eq. (\ref{1}).}

\texttt{Step 4-6- Print $v_{n}(r)$, $\left|v_{n+1}(r)-v_n(r)\right|$
and $\left|v(r) - v_n(r) \right|$. }

\texttt{Step 4-7- $n=n +1$. }

\texttt{ \}}

A sample program of the CADNA library is presented in Appendix
section.

\textbf{Example 1:} Consider the following weakly regular Volterra
IE \cite{dn1}
\begin{equation}\label{6-1}
\begin{array}{l}
\displaystyle  \int_0^{\frac{r}{8}} (1+r+s) v(s) ds +
\int_{\frac{r}{8}}^{\frac{r}{2}} (2+rs) v(s) ds +
\int_{\frac{r}{2}}^{\frac{3r}{4}} (r+s-1) v(s) ds
  - 4 \int_{\frac{3r}{4}}^{r} v(s) ds\\
  \\
\displaystyle  = \frac{1}{128}\Big[ -4 - \frac{1}{8} (16r + 69
r^2+15r^3) - \exp(\frac{r}{4}) (r^2-13r+12)+\exp(r) (4r^2-16r+28) \\
\\
~~~~~~~~~~ \displaystyle + \exp (\frac{3r}{2}) (14r+20) -32 \exp
(2r) \Big],
\end{array}
\end{equation}
where the exact solution is $v(r) = \frac{\exp (2r) -1}{8}$. In
Table \ref{t1}, the numerical results of the CESTAC method for
$r=0.2$ are shown. According to these results, the presented
algorithm stopped when the NCSDs between approximations $v_{n}(r),
v_{n-1}(r)$ are almost equal to NCSDs between $v(r), v_{n}(r)$ and
Theorem \ref{th6} can support the obtained results theoretically.
According to numerical results this equality is obtained in $10$-th
iteration. So the optimal iteration of Taylor-collocation method to
solve the weakly regular Volterra IE of the first kind (\ref{6-1})
is $n=10$. As mentioned before, if we use the Taylor-collocation
method based on the FPA we should apply the stopping condition
(\ref{cond1}) which it depends on parameter $\varepsilon$. The
numerical results of the FPA are presented in table \ref{t3} for
$\varepsilon = 10^{-10}$ and $r=0.2$. Also in table \ref{t2}, the
number of iterations for different values of $\varepsilon$ are
demonstrated. According to this table for large $\varepsilon$, the
iterations are stopped before getting to a proper approximation and
for small $\varepsilon$, the useless iterations are done without
amending the efficiency of the results. Furthermore, in order to
find the numerical approximations based on the FPA, the exact
solution is important because we should apply the termination
criterion (\ref{cond1}) instead of condition (\ref{cond2}). They are
some of most important defects of the FPA than the SA. The sample
program of CADNA library to solve the example 1 is presented in
Appendix.

\begin{table}
\caption{ Numerical results of Example 1 based on the CESTAC method
and the stopping condition (\ref{cond2}) for $r=0.2$. }\label{t1}
 \centering
\scalebox{0.9}{
\begin{tabular}{|c|l|l|l|}
  \hline
$n$   &~~~~~~~~~~$v_n(r)$   & ~~~~~$|v_{n}(r) -v_{n-1}(r) |$  & ~~~~~~~~$|v(r) - v_{n}(r)|$ \\
    \hline
    2    &        0.29415165711104E+001  &    0.29415165711104E+001     &        0.28800384839052E+001\\
    3    &       -0.2374E+000            &    0.317901E+001        &     0.2989E+000\\
    4    &        0.1509E+000            &    0.38844E+000         &    0.8946E-001\\
    5    &        0.4230E-001            &    0.10864E+000         &    0.1917E-001\\
    6    &        0.6480037878687E-001   &             0.2249E-001    &         0.332229158171E-002\\
    7    &        0.6099E-001            &   0.380E-002         &    0.48E-003\\
    8    &        0.6154E-001            &    0.54E-003         &    0.6E-004\\
    9    &        0.6147075816137E-001   &             0.6E-004     &        0.732904378E-005\\
    10   &         0.6147E-001           &     @.0                   &
    @.0\\
 \hline
 \end{tabular}
 }
\end{table}

\begin{table}
\caption{ Numerical results of Example 1 based on the condition
(\ref{cond1}) for $\varepsilon = 10^{-10}$ and $r=0.2$. }\label{t3}
 \centering
\scalebox{0.9}{
\begin{tabular}{|c|l|l|}
  \hline
$n$   & ~~~~~~~~~~$v_n(r)$    & ~~~~~~~$ | v(r) - v_n(r) | $   \\
    \hline
    2    &        2.94151657111040654158     &       2.88003848390524774814\\
    3    &        -0.23749678750420105677    &        0.29897487470935985021\\
    4    &        0.15094246356455431890     &       0.08946437635939552546\\
    5    &        0.04229965285857797269     &       0.01917843434658082075\\
    6    &        0.06480037878687468222     &       0.00332229158171588879\\
    7    &        0.06099350668538974785     &       0.00048458051976904559\\
    8    &        0.06154031969903288324     &       0.00006223249387408980\\
    9    &        0.06147075816137072268     &       0.00000732904378807075\\
    10   &         0.06147890505841218517    &        0.00000081785325339173\\
    11   &         0.06147799914393221182    &        0.00000008806122658161\\
    12   &         0.06147809636759438146    &        0.00000000916243558802\\
    13   &         0.06147808629013713083    &        0.00000000091502166261\\
    14   &         0.06147808728981307008    &        0.00000000008465427664\\
 \hline
 \end{tabular}
 }
\end{table}

\begin{table}
\caption{ Number of iterations of Taylor-collocation method for
Example 1 based on the FPA and termination criterion (\ref{cond1})
for different values of $\varepsilon$ and $r=0.2$. }\label{t2}
 \centering
\scalebox{1}{
\begin{tabular}{|c|c|c|c|c|c|c|c|}
\hline
$\varepsilon$&small values&$10^{-10}$&$10^{-5}$&$10^{-1}$&$0.5$&$1$&large values\\
\hline
$n$&$>>14$&$14$&9&4&3&3&2\\
 \hline
 \end{tabular}
 }
\end{table}

\textbf{Example 2:} Consider the following IE \cite{dn2}
\begin{equation}\label{6-2}
\displaystyle   2 \int_0^{\sin\frac{r}{2}} v(s) ds -
\int_{\sin\frac{r}{2}}^{2 \sin\frac{r}{3}} v(s) ds +  \int_{2
\sin\frac{r}{3}}^{r} v(s) ds
   = \frac{r^3}{3} + \sin^3 \frac{r}{2} - \frac{16}{3}
\sin^3
  \frac{r}{3}, ~~~~~r \in [0,\frac{3 \pi}{2}],
\end{equation}
where $v(r) = r^2$. The numerical results of table \ref{t4} are
obtained based on the SA. As we know, we should apply the stopping
condition (\ref{cond2}). So only two successive approximations are
required and we do not need to have the exact solution. According to
table \ref{t4} the optimal step of presented method is $n=4$ and the
optimal approximation for $r=0.5$  is $v_n(0.5) =
0.24999999999999E+000$. In table \ref{t6}, the numerical results of
the FPA are shown for large $\varepsilon$. Thus its algorithm will
be stopped without getting to suitable approximation with high
accuracy. If we apply the FPA, finding the proper $\varepsilon$ is
important and it is one of faults of computations by FPA. Numerical
results of table \ref{t5} are very good paradigm to show the
mentioned reality.

\begin{table}
\caption{Numerical results of Example 2 by using the CESTAC method
for $r=0.5$.}\label{t4}
 \centering
\scalebox{0.9}{
\begin{tabular}{|c|l|l|l|}
  \hline
$n$   &~~~~~~~~~~$v_n(r)$   & ~~~~~$|v_{n}(r) -v_{n-1}(r) |$  & ~~~~~~~~$|v(r) - v_{n}(r)|$ \\
    \hline
    2    &     0.190380238054409E+001        &         0.190380238054409E+001       &       0.165380238054409E+001\\
    3    &      0.25000000000000E+000         &        0.165380238054409E+001       &       @.0\\
    4    &       0.24999999999999E+000        &         @.0                        &         @.0\\
 \hline
 \end{tabular}
 }
\end{table}

\begin{table}
\caption{Obtained results of Example 2 based on the condition
(\ref{cond1}) for large value of $\varepsilon$ and
$r=0.5$.}\label{t6}
 \centering
\scalebox{0.9}{
\begin{tabular}{|c|l|l|}
  \hline
$n$   & ~~~~~~~~~~$v_n(r)$    & ~~~~~~~$ | v(r) - v_n(r) | $   \\
    \hline
    2    &    1.90380238054409267612   &         1.65380238054409267612\\
 \hline
 \end{tabular}
 }
\end{table}

\begin{table}
\caption{Iterations of Taylor-collocation method for Example 2 by
applying the FPA and criterion (\ref{cond1}) for different values of
$\varepsilon$ and $r=0.5$.}\label{t5}
 \centering
\scalebox{1}{
\begin{tabular}{|c|c|c|c|c|c|c|c|}
\hline $\varepsilon$&small
values&$10^{-10}$&$10^{-5}$&$10^{-1}$&$0.5$&$1$&large values\\
\hline
$n$&$>>3$&3&3&3&3&3&2\\
 \hline
 \end{tabular}
 }
\end{table}

\textbf{Example 3:} Let us consider the following Volterra IE
\begin{equation}\label{6-3}
\int_0^{\frac{r}{4}} (1+r+s) v(s) ds +
\int_{\frac{r}{4}}^{\frac{r}{2}} (2 + rs ) v(s) ds +
\int_{\frac{r}{2}}^{r} (1+r+s) v(s) ds = \frac{31 r^6}{40960} +
\frac{1099r^5}{20480}+\frac{271r^4}{8192}.
\end{equation}
where the exact solution is $v(r) = \frac{r^3}{8}$. The numerical
results of Taylor-collocation method to solve (\ref{6-3}) by using
CESTAC method is presented in table \ref{t7}. We should note that
the final column of table \ref{t7} is only for comparison and we do
not need to find the absolute error in the SA because in Theorem
\ref{th6} we proved $C_{v_n(s),v(s)} = C_{v_n(s),v_{n+1}(s)}.$ Thus
we can apply the termination criterion (\ref{cond2}) which it
depends on two successive approximations. In Table \ref{t9}, the
results of Taylor-collocation method based on the FPA are
demonstrated for $\varepsilon = 10^{-1}$ and $r=0.8$ which we can
not guarantee the obtained results because the numerical results are
without proper precision. Finally, table \ref{t8} are presented to
show the number of iterations of presented method for different for
different values of $\varepsilon$.

\begin{table}
\caption{Numerical results of Example 3 by CESTAC method and
condition (\ref{cond2}) for $r=0.8$.}\label{t7}
 \centering
\scalebox{0.9}{
\begin{tabular}{|c|l|l|l|}
  \hline
$n$   &~~~~~~~~~~$v_n(r)$   & ~~~~~$|v_{n}(r) -v_{n-1}(r) |$  & ~~~~~~~~$|v(r) - v_{n}(r)|$ \\
    \hline
    2    &        0.187158677184107E+000    &            0.187158677184107E+000    &         0.144283677184107E+000\\
    3    &        0.48445393318581E-001     &           0.13871328386552E+000      &       0.5570393318581E-002\\
    4    &        0.42874999999999E-001     &           0.5570393318581E-002       &      @.0\\
    5    &        0.42874999999999E-001     &           0.1E-015      &       @.0\\
    6    &        0.4287499999999E-001      &          @.0            &                    @.0\\
 \hline
 \end{tabular}
 }
\end{table}

\begin{table}
\caption{Results of Example 3 based on the condition (\ref{cond1})
for $\varepsilon = 10^{-1}$ and $r=0.8$.}\label{t9}
 \centering
\scalebox{0.9}{
\begin{tabular}{|c|l|l|}
  \hline
$n$   & ~~~~~~~~~~$v_n(r)$    & ~~~~~~~$ | v(r) - v_n(r) | $   \\
    \hline
        2    &        0.18715867718410730824   &         0.14428367718410731180\\
    3        &    0.04844539331858149778       &     0.00557039331858150827\\
 \hline
 \end{tabular}
 }
\end{table}

\begin{table}
\caption{Iterations of Taylor-collocation method for Example 3 by
using the FPA and criterion (\ref{cond1}) for different values of
$\varepsilon$ and $r=0.8$.}\label{t8}
 \centering
\scalebox{1}{
\begin{tabular}{|c|c|c|c|c|c|c|c|}
\hline $\varepsilon$&small
values&$10^{-10}$&$10^{-5}$&$10^{-1}$&$0.5$&$1$&large values\\
\hline
$n$&$>>4$&$4$&4&3&2&2&2\\
 \hline
 \end{tabular}
 }
\end{table}

\textbf{Example 4:} Abel IEs are one kind of Volterra IEs with
singular kernel which have many applications in science and
engineering \cite{man6,man10,man7,siz5,waz}. In this example, we
consider the following linear Abel IE of the first kind
\begin{equation}\label{6-4}
\int_0^r \frac{v(s)}{\sqrt{r^2-s^2}} ds = \frac{2}{3} \pi r^3,
\end{equation}
where the exact solution is $v(r) = \pi r^3$. According to the
numerical results of Table \ref{t10} the optimal iteration of
Taylor-collocation method is $n=15$ with optimal approximation
$v(0.1)=0.314159E-002$. The numerical results of FPA are obtained
for $\varepsilon = 10^{-5}$ and $r=0.1$ in Table \ref{t11} that
presented algorithm is stopped in 9-th iteration. Also, the number
of iterations for different values of $\varepsilon$ is demonstrated
in Table \ref{t12}.

\begin{table}
\caption{Numerical results of Example 4 by CESTAC method and
condition (\ref{cond2}) for $r=0.1$.}\label{t10}
 \centering
\scalebox{0.9}{
\begin{tabular}{|c|l|l|l|}
  \hline
$n$   &~~~~~~$v_n(r)$   & ~$|v_{n}(r) -v_{n-1}(r) |$  & ~$|v(r) - v_{n}(r)|$ \\
    \hline
    2  &          0.1009201E+001   &             0.1009201E+001  &           0.100605E+001\\
    3  &         -0.368416E+000     &           0.1377617E+001    &         0.371558E+000\\
    4  &          0.2301187E-001    &            0.3914283E+000   &          0.1987028E-001\\
    5  &          0.585469E-002     &           0.171571E-001    &         0.271310E-002\\
    6  &          0.285618E-002     &           0.299850E-002    &         0.28540E-003\\
    7  &          0.3002522E-002    &            0.14633E-003    &         0.13906E-003\\
    8 &           0.3116390E-002     &           0.11386E-003    &         0.2520E-004\\
    9 &           0.313838E-002     &           0.21993E-004    &         0.320E-005\\
    10 &           0.313945E-002     &           0.106E-005     &        0.213E-005\\
    11  &          0.3141559E-002    &            0.296E-006    &         0.31E-007\\
    12  &          0.3141684E-002    &            0.12E-006     &        0.92E-007\\
    13  &          0.314157E-002     &           0.11E-006      &       0.1E-007\\
    14  &          0.314159E-002     &           0.1E-007       &      @.0\\
    15  &          0.314159E-002     &           @.0             &                   @.0\\
 \hline
 \end{tabular}
 }
\end{table}

\begin{table}
\caption{Results of Example 4 based on the condition (\ref{cond1})
for $\varepsilon = 10^{-5}$ and $r=0.1$.}\label{t11}
 \centering
\scalebox{0.9}{
\begin{tabular}{|c|l|l|}
  \hline
$n$   & ~~~~~~~~~~$v_n(r)$    & ~~~~~~~$ | v(r) - v_n(r) | $   \\
    \hline
    2   &         1.00920140743255615234     &       1.00605976581573486328\\
    3   &         -0.36841639876365661621    &        0.37155798077583312988\\
    4   &         0.02301186881959438324     &       0.01987027563154697418\\
    5   &         0.00585469184443354607      &      0.00271309888921678066\\
    6   &         0.00285618798807263374      &      0.00028540496714413166\\
    7   &         0.00300252321176230907     &       0.00013906974345445633\\
    8   &         0.00311639113351702690     &       0.00002520182169973850\\
    9   &         0.00313838454894721508     &       0.00000320840626955032\\
 \hline
 \end{tabular}
 }
\end{table}

\begin{table}
\caption{Iterations of Taylor-collocation method for Example 4 by
using the FPA and criterion (\ref{cond1}) for different values of
$\varepsilon$ and $r=0.1$.}\label{t12}
 \centering
\scalebox{1}{
\begin{tabular}{|c|c|c|c|c|c|c|c|}
\hline $\varepsilon$&small
values&$10^{-10}$&$10^{-5}$&$10^{-1}$&$0.5$&$1$&large values\\
\hline
$n$&$>>9$&$>9$&9&4&3&3&2\\
 \hline
 \end{tabular}
 }
\end{table}

\textbf{Example 5:} Let us consider the following non-linear Abel IE
of the first kind \cite{man6,waz}
\begin{equation}\label{6-5}
\pi+r = \int_0^r \frac{\sin^{-1}(v(s))}{\sqrt{r^2-s^2}} ds,
\end{equation}
where the exact solution is $v(r) = \sin(r+2)$. In order to solve
the mentioned problem at first we can convert the non-linear IE
(\ref{6-5}) to the linear form by using the transformations
\begin{equation}\label{6-6}
v'(s) = \sin^{-1}(v(s)),~~~~~v(s) = \sin (v'(s)).
\end{equation}
Thus, we should solve the following IE
\begin{equation}\label{6-7}
\pi+r = \int_0^r \frac{v'(s)}{\sqrt{r^2-s^2}} ds,
\end{equation}
instead of IE (\ref{6-5}). Then by substituting the obtained
solution in Eq. (\ref{6-6}) we can find the approximate solution of
IE (\ref{6-5}). In Table \ref{t13}, the optimal iteration of
Taylor-collocation method, the optima approximation and the optimal
error are obtained for $r=0.4$. In Table \ref{t14}, because the
results are based on the FPA and in this case we do not know the
optimal value of $\varepsilon$ so the algorithm is stopped in step 6
without getting to suitable approximation. Finally, in Table
\ref{t15}, the number of iterations which are generated by applying
the criterion (\ref{cond1}) are demonstrated for $r=0.4$ and
different values of $\varepsilon$.

\begin{table}
\caption{Numerical results of Example 5 by CESTAC method and
condition (\ref{cond2}) for $r=0.4$.}\label{t13}
 \centering
\scalebox{0.9}{
\begin{tabular}{|c|l|l|l|}
  \hline
$n$   &~~~~~~~$v_n(r)$   & ~$|v_{n}(r) -v_{n-1}(r) |$  & ~$|v(r) - v_{n}(r)|$ \\
    \hline
    2   &        0.7086622E+000    &            0.7086622E+000   &          0.331989E-001\\
    3   &         0.6878617E+000   &             0.208003E-001   &          0.12398E-001\\
    4   &         0.6778912E+000   &             0.99705E-002    &         0.24279E-002\\
    5   &         0.6755348E+000   &             0.23564E-002   &          0.715E-004\\
    6   &         0.6754605E+000   &             0.743E-004     &        0.28E-005\\
    7   &         0.6754629E+000   &             0.24E-005      &       0.3E-006\\
    8   &         0.6754631E+000    &            0.1E-006       &      0.2E-006\\
    9   &         0.6754631E+000    &            @.0            &                    @.0\\
 \hline
 \end{tabular}
 }
\end{table}

\begin{table}
\caption{Results of Example 5 based on the condition (\ref{cond1})
for $\varepsilon = 10^{-5}$ and $r=0.4$.}\label{t14}
 \centering
\scalebox{0.9}{
\begin{tabular}{|c|l|l|}
  \hline
$n$   & ~~~~~~~~~~$v_n(r)$    & ~~~~~~~$ | v(r) - v_n(r) | $   \\
    \hline
    2  &          0.70866230274517871823    &        0.03319919251737557531\\
    3  &          0.68786180640627836436    &        0.01239869617847522143\\
    4  &          0.67789121429022713983    &        0.00242810406242399690\\
    5  &          0.67553501264255699787    &        0.00007190241475385495\\
    6  &          0.67546047309782475399    &        0.00000263712997838894\\
 \hline
 \end{tabular}
 }
\end{table}

\begin{table}
\caption{Iterations of Taylor-collocation method for Example 5 by
using the FPA and criterion (\ref{cond1}) for different values of
$\varepsilon$ and $r=0.4$.}\label{t15}
 \centering
\scalebox{1}{
\begin{tabular}{|c|c|c|c|c|c|c|c|}
\hline $\varepsilon$&small
values&$10^{-10}$&$10^{-5}$&$10^{-1}$&$0.5$&$1$&large values\\
\hline
$n$&$>>6$&$>6$&6&2&2&2&2\\
 \hline
 \end{tabular}
 }
\end{table}

\section{Conclusion}

The SA and specially using the CESTAC method and the CADNA library
are applicable methods and tools to solve the mathematical problems
with controlling the accuracy and detecting the rate of numerical
errors against applying the methods based on the FPA. In this paper,
the Taylor-collocation method to solve the first kind weakly regular
Volterra  IEs with discontinuous kernels was illustrated based on
the SA. The methods based on the FPA have some problems that in
order to void them we should apply the new arithmetic. Thus we
suggested to apply the SA instead of the FPA. In order to find the
numerical results based on the FPA, knowing the exact solution is
important but the SA depends on two successive approximations. The
termination criterion of the FPA depends on small value like
$\varepsilon$ but the SA independents of $\varepsilon$. In the FPA,
we used the mathematical softwares but in the SA, we applied the
CADNA library. By using the CADNA library, not only the optimal
approximation, the optimal iteration and the optimal error can be
obtained but also some of numerical instabilities can be found. Our
license to use the SA instead of the FPA is to prove the main
theorem which in this theorem we showed that the NCSDs between two
approximations $v_{n}(r), v_{n-1}(r)$ are almost equal to the NCSDs
between $v(r), v_n(r)$. Thus, we can apply the SA and the stopping
condition (\ref{cond2}) instead of the FPA and the termination
criterion (\ref{cond1}).


\section*{Appendix}
In this section, a sample CADNA library program to solve the Example
1 is presented. This program is based on the stochastic arithmetic.
The following codes are similar to C++ codes and only some parts are
different. This program is presented is some parts.

\begin{itemize}
  \item In part 1, the commands of C++ are recalled. Specially by using
\texttt{$\sharp$include <cadna.h>} the commands of CADNA library can
be used.
  \item Part 2 until end of program is the main part. We should note that in
order to introduce the variables or arrays the double precision is
applied with $_{-}$st as \texttt{double$_{-}$st} because the
variables will be based on the stochastic arithmetic.
  \item In part 3, we can calculate the nodes of Taylor-collocation method.
In every part, some disable lines are exist that we can show the
results of those part by enabling them.
  \item Part 4 is to construct the array of matrix $A$. In this part, there
are some functions as \texttt{simp3, f, ff} and others that we did
not mention here.
  \item By using part 5, we can find the arrays of right hand matrix $F$.
  \item In part 6, the obtained system of equation is solved.
  \item Part 7 is the final part of this program that we apply to show the
final results. It is important that in order to show the numerical
results with significant digits we should apply \texttt{$\%s$} and
\texttt{strp} in print line. Also command \texttt{while} depends on
difference of two successive approximations. Finally command
\texttt{cadna$_{-}$end()} is the end of CADNA library.
\end{itemize}

\textbf{PART 1:}\\
\texttt{$\sharp$include <stdio.h> \\
$\sharp$include <math.h>\\
$\sharp$include <cadna.h>}

\textbf{PART 2:}\\
\texttt{ main()\\}
\{ \\
\texttt{cadna$_{-}$init(-1);\\
    double$_{-}$st
    z,H[50][50],s[50],a,b,M[50][50],y[50],S[50],exact,T[50],r;\\
int n,i,j,k; \\
a=0;b=1; S[0]=0; \\
n=1; z=0.2;\\
exact=(1./8)*(exp(2*z)-1);\\
printf("-----------------------------------------------------------------
$\backslash$ n"); \\
printf("~~n~~~~~approximate solution~~~~~difference of two
term~~~~~absolute error $\backslash$ n");\\
printf("-----------------------------------------------------------------
$\backslash$ n"); \\
do\\}
    \{\\
\texttt{S[n]=0;}

\textbf{PART 3:}\\
\texttt{for(i=1;i<=n+1;i++)\\}
 \{\\
\texttt{    s[i-1] = a+((b-a)/n)*(i-1);\\} \texttt{//  printf(" \%s
$\backslash$n",strp(s[i-1]));\\} \}

\textbf{PART 4:}\\
\texttt{for(i=1;i<=n+1;i++)\\}
 \{\\
 \texttt{for(j=1;j<=n+1;j++)\\}
  \{\\
\texttt{H[i-1][j-1] =g(s[i-1],j,n);\\} \texttt{//printf("   \%s
$\backslash$ n",strp(H[i-1][j-1]));\\}
\texttt{M[i-1][j-1]=simp3(f,s[i-1],a,s[i-1]*(1./8),j,500)}

 ~~~~~~~~~~~~~\texttt{+simp3(ff,s[i-1],s[i-1]*(1./8),s[i-1]*(1./2),j,500) }

 ~~~~~~~~~~~~~\texttt{+simp3(fff,s[i-1],s[i-1]*(1./2),s[i-1]*(3./4),j,500)}

 ~~~~~~~~~~~~~\texttt{-simp3(ffff,s[i-1],s[i-1]*(3./4),s[i-1],j,500);\\}
\texttt{//printf("   \%s    $\backslash$ n",strp(M[i-1][j-1]));\\}
\}\\
 \}

\textbf{PART 5:}\\
\texttt{for(i=1;i<=n+1;i++)\\}
 \{\\
    \texttt{y[i-1] = (1./128)*(-4-((1./8)*(16*s[i-1]+69*s[i-1]*s[i-1]+15*s[i-1]*s[i-1]*s[i-1]))}

    \texttt{~~~~-(exp(s[i-1]*(1./4))*(s[i-1]*s[i-1]-13*s[i-1]+12))}

    \texttt{~~~~+(exp(s[i-1])*(4*s[i-1]*s[i-1]-16*s[i-1]+28))}

    \texttt{~~~~+(exp((3./2)*s[i-1])*(14*s[i-1]+20))-(32*exp(2*s[i-1])));\\}
\texttt{//  printf("   \%s    $\backslash$ n",strp(y[i-1]));\\} \}

\textbf{PART 6:}\\
\texttt{double$_{-}$st temp;
\\} \texttt{for(k=0;k<=n;k++)\\}
 \{\\
\texttt{ temp = M[k][k]; \\}
 \texttt{M[k][k] = 1.0;\\}
\texttt{for(j=0;j<=n;j++) \\} \texttt{M[k][j]=M[k][j]/temp; \\}
\texttt{for(i=0;i<=n;i++)\\}
 \texttt{if(i!=k)\\}
  \{\\
 \texttt{ temp=M[i][k]; \\}
  \texttt{M[i][k] = 0.0;\\}
  \texttt{for(j=0;j<=n;j++) \\}
  \texttt{M[i][j] =M[i][j]- temp*M[k][j];\\}
   \}\\
    \}\\
\texttt{//for(i=0;i<=n;i++)\\} \texttt{//\{for(j=0;j<=n;j++)\\}
\texttt{//printf("    \%s     ",strp(M[i][j]));\\}
\texttt{//printf("$\backslash$ n");\}\\}
\texttt{double$_{-}$st e,c[50];\\}
 \texttt{for(i=0;i<=n;i++)\\}
  \{\\
  \texttt{ e=0;\\}
\texttt{for(j=0;j<=n;j++)\\}
 \{\\
  \texttt{e=e+M[i][j]*y[j];\\}
   \}\\
   \texttt{ c[i]=e;\\}
\texttt{//printf("$\backslash$ n   \%.10f   $\backslash$n",c[i]);\\}
\}

\textbf{PART 7:}\\
\texttt{for(i=0;i<=n;i++)\\}
 \{\\
 \texttt{ S[n]=S[n]+c[i]*g(z,i+1,n);\\}
   \}\\
\texttt{printf("    \%d           \%s               \%s \%s
$\backslash$n",n+1,strp(S[n]),strp(fabs(S[n]-S[n-1]))}

\texttt{~~~~~,strp(fabs(exact-S[n])));\\} \texttt{n++;\\}
\}\\
\texttt{ while(S[n-1]-S[n-2]!=0);\\}
\texttt{printf("----------------------------------------------------------$\backslash$
n $\backslash$ n $\backslash$ n");\\}
\texttt{cadna$_{-}$end();\\}
 \}\\
\end{document}